\documentclass[12pt]{article}
\usepackage{amsmath, amsthm, amsfonts}
\usepackage{amssymb}
\usepackage{amscd}
\usepackage{verbatim}
\begin{document}
\newcommand{\loc}{{\mathrm{loc}}}
\newcommand{\dx}{\,\mathrm{d}x}
\newcommand{\dy}{\,\mathrm{d}y}
\newcommand{\dz}{\,\mathrm{d}z}
\newcommand{\dt}{\,\mathrm{d}t}
\newcommand{\du}{\,\mathrm{d}u}
\newcommand{\dv}{\,\mathrm{d}v}
\newcommand{\dV}{\,\mathrm{d}V}
\newcommand{\ds}{\,\mathrm{d}s}
\newcommand{\dr}{\,\mathrm{d}r}
\newcommand{\dS}{\,\mathrm{d}S}
\newcommand{\drho}{\,\mathrm{d}\rho}
\newcommand{\core}{C_0^{\infty}(\Omega)}
\newcommand{\sob}{W^{1,p}(\Omega)}
\newcommand{\sobloc}{W^{1,p}_{\mathrm{loc}}(\Omega)}
\newcommand{\merhav}{{\mathcal D}^{1,p}}
\newcommand{\be}{\begin{equation}}
\newcommand{\ee}{\end{equation}}
\newcommand{\mysection}[1]{\section{#1}\setcounter{equation}{0}}
\newcommand{\bea}{\begin{eqnarray}}
\newcommand{\eea}{\end{eqnarray}}
\newcommand{\bean}{\begin{eqnarray*}}
\newcommand{\eean}{\end{eqnarray*}}
\newcommand{\thkl}{\rule[-.5mm]{.3mm}{3mm}}
\newcommand{\cw}{\stackrel{\rightharpoonup}{\rightharpoonup}}
\newcommand{\id}{\operatorname{id}}
\newcommand{\supp}{\operatorname{supp}}
\newcommand{\wlim}{\mbox{ w-lim }}
\newcommand{\mymu}{{x_N^{-p_*}}}
\newcommand{\R}{{\mathbb R}}
\newcommand{\N}{{\mathbb N}}
\newcommand{\Z}{{\mathbb Z}}
\newcommand{\Q}{{\mathbb Q}}
\newcommand{\abs}[1]{\lvert#1\rvert}
\newtheorem{theorem}{Theorem}[section]
\newtheorem{corollary}[theorem]{Corollary}
\newtheorem{lemma}[theorem]{Lemma}
\newtheorem{definition}[theorem]{Definition}
\newtheorem{remark}[theorem]{Remark}
\newtheorem{proposition}[theorem]{Proposition}
\newtheorem{problem}[theorem]{Problem}
\newtheorem{conjecture}[theorem]{Conjecture}
\newtheorem{question}[theorem]{Question}
\newtheorem{example}[theorem]{Example}
\newtheorem{Thm}[theorem]{Theorem}
\newtheorem{Lem}[theorem]{Lemma}
\newtheorem{Pro}[theorem]{Proposition}
\newtheorem{Def}[theorem]{Definition}
\newtheorem{Exa}[theorem]{Example}
\newtheorem{Exs}[theorem]{Examples}
\newtheorem{Rems}[theorem]{Remarks}
\newtheorem{Rem}[theorem]{Remark}
\newtheorem{Cor}[theorem]{Corollary}
\newtheorem{Conj}[theorem]{Conjecture}
\newtheorem{Prob}[theorem]{Problem}
\newtheorem{Ques}[theorem]{Question}
\newcommand{\pf}{\noindent \mbox{{\bf Proof}: }}
\renewcommand{\theequation}{\thesection.\arabic{equation}}
\catcode`@=11
\@addtoreset{equation}{section}
\catcode`@=12
\title{A Liouville-type theorem for the $p$-Laplacian with potential term}
\author{Yehuda Pinchover\\
 {\small Department of Mathematics}\\ {\small  Technion - Israel Institute of Technology}\\
 {\small Haifa 32000, Israel}\\
{\small pincho@techunix.technion.ac.il}\and  Achilles Tertikas
\\{\small  Department
of Mathematics}\\{\small University of Crete}\\{\small 714 09
Heraklion, GREECE}\\{\small tertikas@math.uoc.gr}
 \and  Kyril Tintarev
\\{\small Department of Mathematics}\\{\small Uppsala University}\\
{\small SE-751 06 Uppsala, Sweden}\\{\small
kyril.tintarev@math.uu.se}}
\maketitle
\newcommand{\dnorm}[1]{\thkl #1 \thkl\,}
\thispagestyle{empty}
\begin{abstract} In this paper we prove a sufficient condition, in terms of
the behavior of a ground state of a singular $p$-Laplacian problem
with a potential term, such that a nonzero subsolution of another
such problem is also a ground state. Unlike in the linear case
($p=2$), this condition involves comparison of both the functions
and of their gradients.
\\[1mm]
\end{abstract}
 \mysection{Introduction}\label{secint}
Positivity properties of quasilinear elliptic equations, in
particular those with the $p$-Laplacian term in the principal
part, have been extensively studied over the recent decades (see
for example \cite{AH1,AH2,HKM,ky3} and the references therein).
Fix $p\in(1,\infty)$, and a domain $\Omega\subseteq\R^d$.  In this
paper we use positivity properties of such equations to prove a
general Liouville comparison principle for equations of the form
$$-\Delta_p(u)+V|u|^{p-2}u=0\quad \mbox{in }  \Omega,$$
where  $\Delta_p(u):=\nabla\cdot(|\nabla u|^{p-2}\nabla u)$ is the
$p$-Laplacian, and $V\in L_\mathrm{loc}^\infty(\Omega;\R)$ is a
given potential. Throughout this paper we assume that \be
\label{Q} Q(u):=\int_\Omega \left(|\nabla u|^p+V|u|^p\right)\dx\ge
0\ee for all $u\in \core$.
\begin{definition}{\em We say that a function
$v\in W^{1,p}_{\mathrm{loc}}(\Omega)$ is a {\em (weak) solution}
of the equation  \be \label{groundstate}
 \frac{1}{p}Q^\prime
(v):=-\Delta_p(v)+V|v|^{p-2}v=0\quad \mbox{in }  \Omega,\ee if for
every $\varphi\in\core$
 \be \label{solution} \int_\Omega (|\nabla v|^{p-2}\nabla
v\cdot\nabla\varphi+V|v|^{p-2}v\varphi)\dx=0. \ee
We say that a real function $v\in C^1_{\mathrm{loc}}(\Omega)$ is a
{\em supersolution} (resp. {\em subsolution})  of the equation
(\ref{groundstate}) if for every nonnegative $\varphi\in\core$
 \be\label{supersolution}
\int_\Omega (|\nabla v|^{p-2}\nabla
v\cdot\nabla\varphi+V|v|^{p-2}v\varphi)\dx\geq 0 \mbox{ (resp.
}\leq 0\mbox{).} \ee
 }\end{definition}
\begin{remark}{\em It is well-known that any weak solution of (\ref{groundstate})
admits H\"older continuous first derivatives, and that any
nonnegative solution of (\ref{groundstate}) satisfies the Harnack
inequality \cite{Serrin1,Serrin2,T}.
 }\end{remark}
\begin{definition}{\em We say that the functional $Q$ has a
{\em weighted spectral gap in $\Omega$}  if there is a positive
continuous function $W$ in $\Omega$ such that \be\label{wsg}
Q(u)\ge \int_\Omega W|u|^p\dx \qquad \forall u\in\core.\ee
}\end{definition}
\begin{definition}{\em  Let $Q$ be a nonnegative functional on
$\core$. We say that a sequence $\{u_k\}\subset\core$  of
nonnegative functions is a {\em null sequence} of the functional
$Q$ in $\Omega$, if there exists an open set $B\Subset\Omega$
(i.e., $\overline{B}$ is compact in $\Omega$) such that
$\int_B|u_k|^p\dx=1$, and \be
\lim_{k\to\infty}Q(u_k)=\lim_{k\to\infty}\int_\Omega (|\nabla
u_k|^p+V|u_k|^p)\dx=0.\ee
 We say that a positive function $v\in
C^1_{\mathrm{loc}}(\Omega)$ is a {\em ground state} of the
functional $Q$ in $\Omega$ if $v$ is an
$L^p_{\mathrm{loc}}(\Omega)$ limit of a null sequence of $Q$.
}
\end{definition}
\begin{remark}\label{remc1}{\em
The requirement that $\{u_k\}\subset \core$, can clearly be
weakened by assuming only that $\{u_k\}\subset W^{1,p}_0(\Omega)$.
Also, the requirement that $\int_B|u_k|^p\dx=1$ can be replaced by
$\int_B|u_k|^p\dx\asymp 1$, where $f_k\asymp g_k$ means that there
exists a positive constant $C$ such that $C^{-1}g_k\leq f_k \leq
Cg_k$  for all $k\in \mathbb{N}$.
 }\end{remark}
%
The following theorem was proved in \cite{ky3}.
\begin{theorem}
\label{ky3} Let $\Omega\subseteq\R^d$ be a domain, $V\in
L_\mathrm{loc}^\infty(\Omega)$, and $p\in(1,\infty)$. Suppose that
the functional $Q$ is nonnegative on $\core$. Then
\begin{itemize}
\item[(a)] $Q$ has either a weighted spectral gap or a ground
state.

\item[(b)] If $Q$ admits a ground state $v$, then $v>0$ and $v$ satisfies
(\ref{groundstate}).

\item[(c)]
The functional $Q$ admits a ground state if and only if
(\ref{groundstate}) admits a unique positive supersolution.
\end{itemize}
\end{theorem}
\begin{example}\label{ex1}{\em Consider the functional
$Q(u):=\int_{\mathbb{R}^d}|\nabla u|^p\dx$. It follows from
\cite[Theorem~2]{MP} that if $d\leq p$, then $Q$ admits a ground
state $\varphi=\mathrm{constant}$ in $\mathbb{R}^d$. On the other
hand, if $d>p$, then
$$u(x):=\left[1+|x|^{p/(p-1)}\right]^{(p-d)/p}, \qquad v(x)=\mathrm{constant}$$
are two positive supersolutions of the equation $-\Delta_pu=0$ in
$\mathbb{R}^d$. Therefore, Theorem~\ref{ky3} (c) implies that if
$d>p$, then $Q$ has a weighted spectral gap in $\mathbb{R}^d$. See
also Example~\ref{ex2}.
 }\end{example}
In a recent paper \cite{Pl}, Theorem~\ref{ky3} was used in order
to prove, for $p=2$, the following Liouville-type statement.
\begin{theorem}[\cite{Pl}]\label{thm:p=2} Let $\Omega$ be a domain in $\R^d$, $d\geq 1$.
Consider two strictly elliptic Schr\"odinger operators defined on
$\Omega$ of the form
\begin{equation}\label{eqpj}
P_j:=-\nabla\cdot(A_j\nabla)+V_j\qquad j=0,1,
\end{equation}
where  $V_j\in L^{p}_{\mathrm{loc}}(\Omega;\mathbb{R})$ for some
$p>{d}/{2}$, and $A_j:\Omega \rightarrow \mathbb{R}^{d^2}$ are
measurable symmetric matrices such that for any $K\Subset\Omega$
there exists $\mu_K>1$ such that \be \label{stell}
\mu_K^{-1}I_d\le A_j(x)\le \mu_K I_d \qquad \forall x\in K.
 \ee
(Here $I_d$ is the $d$-dimensional identity matrix, and the matrix
inequality $A\leq B$ means that $B-A$ is a nonnegative matrix on
$\mathbb{R}^d$.)

\vskip 3mm

 \noindent Assume that the following assumptions hold true.
\begin{itemize}
\item[(i)] The operator  $P_1$ admits a ground state
 $\varphi$ in $\Omega$.

\item[(ii)]  $P_0\geq 0$ on $\core$, and there exists a
real function $\psi\in H^1_{\mathrm{loc}}(\Omega)$ such that
$\psi_+\neq 0$, and $P_0\psi \leq 0$ in $\Omega$, where
$u_+(x):=\max\{0, u(x)\}$.

\item[(iii)] The following matrix inequality holds
\begin{equation}\label{psialephia}
(\psi_+)^2(x) A_0(x)\leq C\varphi^2(x) A_1(x)\qquad  \mbox{ a. e.
in } \Omega,
\end{equation}
where $C>0$ is a positive constant.
\end{itemize}
Then the operator $P_0$ admits a ground state in $\Omega$, and
$\psi$ is the corresponding ground state. In particular, $\psi$
is (up to a multiplicative constant) the unique positive
supersolution of the equation $P_0u=0$ in $\Omega$.
\end{theorem}
The purpose of this paper is to find an analog of
Theorem~\ref{thm:p=2} when $p\neq 2$. The main statement is as
follows.
\begin{theorem}\label{thm:main} Let $\Omega$ be a domain in $\R^d$, $d\geq 1$, and let
$p\in(1,\infty)$. For $j=0,1$, let $V_j\in L^\infty_\loc(\Omega)$,
and  let
$$ Q_j(u):=\int_\Omega\left(|\nabla
u(x)|^p+V_j(x)|u(x)|^p\right)\dx \qquad u\in\core. $$
\par
 Assume that the following assumptions hold true.
\begin{itemize}
\item[(i)] The functional  $Q_1$ admits a ground state
 $\varphi$ in $\Omega$.

\item[(ii)]  $Q_0\geq 0$ on $\core$, and the equation $Q_0'(u)=0$ in $\Omega$ admits a
subsolution $\psi\in W^{1,p}_{\mathrm{loc}}(\Omega)$ satisfying
$\psi_+\neq 0$.

\item[(iii)] The following inequality holds almost everywhere in $\Omega$
\begin{equation} \label{ineq:phi-psi} \psi_+\le C\varphi, \end{equation} where $C>0$ is a
positive constant.

\item[(iv)] The following inequality holds almost everywhere  in $\Omega\cap\{\psi>0\}$
\begin{equation} \label{ineq:gradphi-gradpsi}
    |\nabla\psi|^{p-2}\le C|\nabla
\varphi|^{p-2},
 \end{equation} where $C>0$ is a positive constant.
\end{itemize}
Then the functional $Q_0$ admits a ground state in $\Omega$, and
$\psi$ is the ground state. In particular, $\psi$ is (up to a
multiplicative constant) the unique positive supersolution of the
equation $Q'_0(u)=0$ in $\Omega$.
\end{theorem}


%
\begin{remark}{\em
Condition \eqref{ineq:gradphi-gradpsi} is redundant for $p=2$. For
$p\neq 2$ it is equivalent to the assumption that the following
inequality holds in $\Omega$:
\begin{equation} \label{ineq:gradphi-gradpsi1}
 \begin{cases}
    |\nabla\psi_+|\le C|\nabla
\varphi| & \text{ if }\; p>2, \\
    |\nabla\psi_+|\ge C|\nabla
\varphi| & \text{ if }\; p<2,
  \end{cases}
 \end{equation} where $C>0$ is a positive constant.}
\end{remark}
\begin{remark}\label{remgradcond}{\em
This theorem holds if, in addition to  \eqref{ineq:phi-psi}, one
assumes instead of $|\nabla\psi|^{p-2}\le C|\nabla \varphi|^{p-2}$
in $\Omega\cap\{\psi>0\}$ (see \eqref{ineq:gradphi-gradpsi}), that
the following inequality holds true almost everywhere in
$\Omega\cap\{\psi>0\}$
\begin{equation} \label{ineq:gradphi-gradpsi-a}
    \psi^2|\nabla\psi|^{p-2}\le C\varphi^2|\nabla
\varphi|^{p-2},
 \end{equation}
 where $C>0$ is a positive constant.
This can be easily observed by repeating the proof of Theorem~\ref{thm:main}
with the equivalent energy functional represented in the form \eqref{allp} instead of \eqref{p<2}.
 }
 \end{remark}
\begin{remark}\label{remNogradcond} {\em Suppose that $1<p<2$, and assume that the ground state
$\varphi>0$ of the functional $Q_1$ is such that $w=\mathbf{1}$ is
a ground state of the functional \be \label{E-1}
E_1^{\varphi}(w)=\int_\Omega \varphi^p|\nabla w|^p \dx , \ee that
is, there is a sequence $\{w_k\}\subset C_0^\infty(\Omega)$ of
nonnegative functions satisfying $E_1^{\varphi}(w_k)\to 0$, and
$\int_B |w_k|^p=1$ for a fixed $B\Subset\Omega$ (this implies that
$w_k\to\mathbf{1}$ in $L^p_{\mathrm{loc}}(\Omega)$). In this case,
the conclusion of Theorem~\ref{thm:main} holds if there is a
nonnegative subsolution $\psi_{+}$ of $Q'_0(u)=0$ satisfying
\eqref{ineq:phi-psi} alone, without any assumption on the
gradients (like \eqref{ineq:gradphi-gradpsi} or
\eqref{ineq:gradphi-gradpsi-a}). This statement follows from the
proof of Theorem~\ref{thm:main} together with the trivial
inequality
$$\int_\Omega v^2 |\nabla
w|^2\left(w|\nabla v|+v|\nabla w|\right)^{p-2}\dx\leq \int_\Omega
v^p|\nabla w|^p \dx$$ which actually holds pointwise. We use this
observation in Example~\ref{ex2}.}
\end{remark}
\begin{remark}{\em
 By Picone identity, a nonnegative functional $Q$ can be
represented as the integral of a nonnegative Lagrangian $L$.
Although the expression for $L$ contains  an indefinite term (see
\eqref{piconeLag}), it admits a two-sided estimate by a simplified
Lagrangian  with nonnegative terms (see
Lemma~\ref{lem:superPicone}). We call the functional associated
with this simplified Lagrangian the {\em simplified energy}. It
plays a crucial role in the proof of Theorem~\ref{thm:main}.}
 \end{remark}
\begin{remark}{\em
Condition \eqref{ineq:gradphi-gradpsi} is  essential when $p>2$,
and presumably also when $p<2$. When $p>2$, $\Omega=\R^d$ and $V$
is radially symmetric, Proposition~\ref{prop:irreducible} shows
that the simplified energy functional is not equivalent to either
of its two terms that lead to conditions \eqref{ineq:phi-psi} and
\eqref{ineq:gradphi-gradpsi}, respectively (see also
Remark~\ref{rem:irreducible}).
 }
 \end{remark}
The outline of the paper is as follows. In Section~\ref{secpicone}
we study the representation of $Q$ as a functional with a positive
Lagrangian, and derive the equivalent simplified energy.
Theorem~\ref{thm:main} is proved in Section~\ref{secmain}, and
Section~\ref{secsimpl} is devoted to the irreducibility of the
simplified energy to either of its terms. In
Section~\ref{secappl}, we study a connection between the ground
states of the functional $Q$ and of its linearization.

\mysection{Picone identity}\label{secpicone}
Let $v>0$, $v\in C^1_{\mathrm{loc}}(\Omega)$, and $u\geq 0$,
$u\in\core$. Denote \be R(u,v):=|\nabla u|^p- \nabla
\left(\frac{u^{p}}{v^{p-1}}\right)\cdot|\nabla v|^{p-2}\nabla v,
\ee and
 \be\label{piconeLag} L(u,v):=|\nabla u|^p+(p-1)\frac{u^p}{v^p}|\nabla
v|^p-p\frac{u^{p-1}}{v^{p-1}}\nabla u\cdot|\nabla v|^{p-2}\nabla
v.\ee Then the following {\em (generalized) Picone identity} holds
\cite{DS,AH1,AH2} \be\label{picone0} R(u,v)=L(u,v). \ee
 Write
$L(u,v)=L_1(u,v)+L_2(u,v)$, where
 \be \label{L1}
L_1(u,v):=|\nabla u|^p+(p-1)\frac{u^p}{v^p}|\nabla
v|^p-p\frac{u^{p-1}}{v^{p-1}}|\nabla u||\nabla v|^{p-1},\ee
 and
 \be\label{L2}
 L_2(u,v):=p\frac{u^{p-1}}{v^{p-1}}|\nabla
v|^{p-2}(|\nabla u||\nabla v|-\nabla u\cdot\nabla v)\ge 0.\ee
 From the obvious inequality $t^p+(p-1)-pt\ge 0$, we also have that
$L_1(u,v)\ge 0$. Therefore, $L(u,v)\ge 0$ in $\Omega$. Let $v\in
C_{\mathrm{loc}}^1(\Omega)$ be a positive solution (resp.
subsolution) of (\ref{groundstate}). Using (\ref{picone0}) and
(\ref{solution}) (resp. (\ref{supersolution})),
 we infer that for every $u\in\core$, $u\ge 0$,
 \be\label{QL} Q(u)=\int_\Omega L(u,v)\dx,\qquad \mbox{resp. }\; Q(u) \leq \int_\Omega L(u,v)\dx.
 \ee
Let now $w:=u/v$, where $v$ is a positive solution  of
(\ref{groundstate}) and $u\in\core$, $u\ge 0$. Then \eqref{QL}
implies
\be \label{QL1} Q(vw)=\int_\Omega\left[|v\nabla w+w\nabla v
|^p-w^p|\nabla v|^p -pw^{p-1}v|\nabla v|^{p-2}\nabla v\cdot\nabla
w\right]\dx. \ee Similarly,  if $v$ is a nonnegative subsolution
of (\ref{groundstate}), then \be \label{QL1sub} Q(vw)\leq
\int_\Omega\left[|v\nabla w+w\nabla v |^p-w^p|\nabla v|^p
-pw^{p-1}v|\nabla v|^{p-2}\nabla v\cdot\nabla w\right]\dx. \ee

 A need to study the linearized operator arises
at a certain step in this paper. This linearized operator is a
Schr\"odinger operator of the form \be \label{divform}
Pu:=(-\nabla\cdot(A\nabla)+V)u\qquad \mbox{ in } \Omega. \ee
  We assume that $V\in L^{\infty}_{\mathrm{loc}}(\Omega;\mathbb{R})$,
 and $A:\Omega \rightarrow \mathbb{R}^{d^2}$ is a
measurable (symmetric) matrix valued function satisfying
\eqref{stell}.   We consider the quadratic form \be \label{assume}
\mathbf{a}[u]:=\int_\Omega\left(A\nabla u\cdot \nabla
u+V|u|^2\right)\mathrm{d}x   \ee on $\core$
 associated with the operator $P$.
  We have the
 following version of Picone identity (see \cite{ky3}).
\begin{lemma}\label{lemqf}  Let $\psi$ be a (real valued) solution of the equation
$P\psi=0$ in $\Omega$. Then for any $v\in \core$ we have
\begin{equation}
\label{e2}\mathbf{a}[\psi v]=
 \int_\Omega \psi^2
A\nabla v\cdot\nabla v \,\mathrm{d}x.
 \end{equation}
 Moreover, if $\psi\in H^1_{\mathrm{loc}}(\Omega)$ is
a nonnegative subsolution of the equation $P\psi=0$ in $\Omega$,
then for any nonnegative $v\in \core$ we have
 \begin{equation}
\label{e1}\mathbf{a}[\psi v]\leq
 \int_\Omega \psi^2
A\nabla v\cdot\nabla v \,\mathrm{d}x.
 \end{equation}
  \end{lemma}
So, in the linear case, the quadratic form induces a convenient
weighted (Dirichlet-type) norm $\|v\|^2:=\int_\Omega \psi^2
A\nabla v\cdot\nabla v \,\mathrm{d}x$ on $\core$.  Recall that in
the quasilinear case ($p\neq 2$), the Lagrangian $L$ in Picone's
identity and \eqref{QL1} contain indefinite terms. Therefore, it
is more convenient to replace identity \eqref{QL1} by two-sided
inequalities with a simpler expression which we call the {\em
simplified  energy}.
\begin{lemma}
\label{lem:superPicone}  Let $v\in C_{\mathrm{loc}}^1(\Omega)$ be
a positive solution of (\ref{groundstate}) and let $w\in
C^1_0(\Omega)$ be a nonnegative function. Then
\bea \label{p<2}  Q(vw)\asymp  \int_\Omega v^2 |\nabla
w|^2\left(w|\nabla v|+v|\nabla w|\right)^{p-2}\dx.\eea

Moreover, for all $p\neq 2$
\bea \label{allp}  Q(vw)\asymp  C\int_{\Omega} |\nabla
w|^2\left(w|\nabla v|v^\frac{2}{p-2}+v^\frac{p}{p-2}|\nabla
w|\right)^{p-2}\dx\eea In particular, for $p>2$ we have
\bea \label{p>2}  Q(vw)\asymp  \int_\Omega  \left(v^p|\nabla
w|^p+v^2|\nabla v|^{p-2} w^{p-2}|\nabla w|^2\right)\dx. \eea

If $v$ is only a nonnegative subsolution of (\ref{groundstate}),
then
\bea \label{p<2sub}  Q(vw)\leq C  \int_{\Omega\cap\{v> 0\}} v^2 |\nabla
w|^2\left(w|\nabla v|+v|\nabla w|\right)^{p-2} \dx.\eea

If $p\neq 2$, then \bea \label{allpsub}  Q(vw)\leq
C\int_{\Omega\cap\{v> 0\}} |\nabla w|^2\left(w|\nabla
v|v^\frac{2}{p-2}+v^\frac{p}{p-2}|\nabla w|\right)^{p-2}\dx\eea
Moreover, for $p>2$ we have
\bea \label{p>2sub}  Q(vw)\leq C  \int_\Omega \left(v^p|\nabla
w|^p+v^2|\nabla v|^{p-2} w^{p-2}|\nabla w|^2\right) \dx. \eea
\end{lemma}
\begin{proof}
Let $1<p<\infty$. We need the following elementary algebraic
vector inequality (cf. \cite{BFT,Shafrir})
 \be \label{ineq p 7}
 |a+b|^p-|a|^p-p|a|^{p-2}a\cdot
b\asymp  |b|^2(|a|+|b|)^{p-2} \ee for all $a,b\in\R^d$.

Indeed, let $t=|b|/|a|$ and $\theta=(a\cdot b)/(|a||b|)$. Note
that for $-1\leq \theta\leq 1$ \be\label{near infty}
\lim_{t\to\infty}\frac{|t^2+2\,\theta t+1|^{p/2}-1-p\,\theta
t}{t^2(1+t)^{p-2}}= 1,\ee and
  \be\label{near zero} \lim_{t\to
0_+}\frac{|t^2+2\,\theta t+1|^{p/2}-1-p\,\theta
t}{t^2(1+t)^{p-2}}= \frac{p}{2}(1+(p-2)\theta^2)>C_p> 0.\ee
Finally, we claim that for $t>0$ and $-1\leq \theta\leq 1$ we have
\be \label{in between} f(t,\theta):=|t^2+2\,\theta
t+1|^{p/2}-1-p\,\theta t> 0. \ee Indeed, set $s:=(t^2+2\,\theta
t+1)^{1/2}\geq 0$, then
$$f(t,\theta)= [s^p+(p-1)-ps] + p[(t^2+2\,\theta t+1)^{1/2}-1-\theta \,
t].$$ Clearly, for $s\geq 0$ we have, $g(s):=[s^p+(p-1)-ps]\geq
0$, and $g(s)=0$ if and only if $s=1$, which holds if and only if
$t=-2\theta$.

On the other hand, let
$$h(t,\theta):=p[(t^2+2\,\theta
t+1)^{1/2}-1-\theta \, t].$$ Then $h(t,\theta)\geq 0$, and
$h(t,\theta)=0$ if and only if $\theta =\pm 1$. Note that if
$\theta =-1$ and $t=-2\theta$, then we have $f(2,-1)= 2p>0$.
 Thus, $f(t,\theta)> 0$ for all $t>0$ and $-1\leq
\theta\leq 1$.

Therefore, for $1<p<\infty $, relations
 \eqref{near infty}--\eqref{in between} imply
$$|t^2+2\,\theta t+1|^{p/2}-1-p\,\theta t\asymp t^2(1+t)^{p-2}.
$$
Thus, \eqref{ineq p 7} holds true for all $a,b\in\R^d$.

Set now $a:=w|\nabla v|$, $b:=v|\nabla w|$. Then we obtain
\eqref{p<2} and \eqref{p<2sub} by applying
 \eqref{ineq p 7} to \eqref{QL1} and \eqref{QL1sub}, respectively.
\end{proof}
The following Allegretto-Piepenbrink-type theorem was proved in
\cite{ky3}.
\begin{theorem}[{\cite[Theorem~2.3]{ky3}}]\label{pos} Let Q be a functional of the form (\ref{Q}).
Then the following assertions are equivalent
\begin{itemize}
\item[(i)] The functional $Q$ is nonnegative on $C_0^\infty(\Omega)$.

\item[(ii)] Equation (\ref{groundstate}) admits a global positive
solution. \item[(iii)] Equation (\ref{groundstate}) admits a
global positive supersolution.
\end{itemize}
\end{theorem}
The next lemma is well known for $p=2$ (see for example
\cite[Lemma~2.9]{Ag}).
\begin{lemma}\label{subsol}
Let $v\in C^1_{\mathrm{loc}}(\Omega)$ be a subsolution  of
equation (\ref{groundstate}). Then $v_+$ is also a subsolution of
(\ref{groundstate}).
\end{lemma}
\begin{proof} Fix $\varphi \in
\core$, $\varphi\geq 0$. As in \cite[Lemma~2.9]{Ag}, define for
$\varepsilon>0$
$$v_\varepsilon:=(v^2+\varepsilon^2)^{1/2}, \qquad  \mbox{ and }\;
\varphi_\varepsilon:=\frac{v_\varepsilon+v}{2v_\varepsilon}\varphi\,.$$
Then $v_\varepsilon \to |v|$, $\nabla v_\varepsilon \to \nabla
|v|$, and $\varphi_\varepsilon \to (\mathrm{sgn}\, v_+) \varphi$
as $\varepsilon \to 0$. An elementary computation shows that
$$\nabla v_\varepsilon\cdot \nabla \varphi \leq \nabla v\cdot
\nabla\left (\frac{v\varphi}{v_\varepsilon}\right),$$ and
therefore,
 $$\nabla
 \left(\frac{v_\varepsilon +v}{2}\right)\cdot \nabla \varphi \leq \nabla
v\cdot  \nabla \varphi_\varepsilon.$$ Since
\be\label{subsolution1} \int_\Omega (|\nabla v|^{p-2}\nabla
v\cdot\nabla\varphi_\varepsilon+V|v|^{p-2}v\varphi_\varepsilon)\dx\leq
0,  \ee it follows that \be\label{subsolution2} \int_\Omega
(|\nabla v|^{p-2}\nabla \left(\frac{v_\varepsilon
+v}{2}\right)\cdot \nabla
\varphi+V|v|^{p-2}v\varphi_\varepsilon)\dx\leq 0.  \ee Letting
$\varepsilon \to 0$ we obtain \be\label{supersolution3}
\int_\Omega (|\nabla v_+|^{p-2}\nabla
v_+\cdot\nabla\varphi+V|v_+|^{p-2}v_+\varphi)\dx\leq 0. \ee
\end{proof}
\mysection{Proof of the main result}\label{secmain}
\begin{proof}[{Proof of Theorem~\ref{thm:main}}]
 By Lemma \ref{subsol}, we may assume that $\psi \geq 0$.

 Let $\{u_k\}$ be a null sequence for $Q_1$, that is
$Q_1(u_k)\to 0$ and, for some nonempty open set $B\Subset\Omega$,
$\int_B u_k^p\dx=1$. Without loss of generality, we may assume
that $B\subset\supp\psi$. Let $w_k:=u_k/\varphi$. From \eqref{p<2}
it follows that with some $C>0$
$$
\int_\Omega \varphi^2 |\nabla w_k|^2\left(w_k|\nabla
\varphi|+\varphi|\nabla w_k|\right)^{p-2}\dx \le C Q_1(u_k)\to 0.
$$
Fix $\alpha,\beta\in \R_+$, then the function $f:\R_+^2\to \R_+$
defined by
$$f(s,t):=\alpha^2 t^2(\beta s^{1/(p-2)}+\alpha t)^{p-2}$$ is nondecreasing monotone function in
each variable separately. Hence, assumptions \eqref{ineq:phi-psi}
and \eqref{ineq:gradphi-gradpsi} imply that
$$
\int_\Omega \psi^2 |\nabla w_k|^2\left(w_k|\nabla
\psi|+\psi|\nabla w_k|\right)^{p-2}\dx \to 0.
$$
Together with \eqref{p<2} this implies that $Q_0(\psi w_k)\to 0$.
On the other hand,  since $w_k\to 1$ in $L^p_\loc(\Omega)$, it
follows that $\psi w_k\to \psi$ in $L^p_\loc(\Omega)$.
Consequently,  $\int_B \varphi^pw_k^p\dx=1$ implies that $\int_B
\psi^p w_k^p\dx\asymp 1$.  In light of Remark~\ref{remc1}, we
conclude that $\psi$ is a ground state of $Q_0$.
\end{proof}

\begin{example}\label{ex12}{\em
Assume that $1\leq d\leq p\leq 2$, $p>1$,   $\Omega
=\mathbb{R}^d$, and consider the functional
$Q_1(u):=\int_{\mathbb{R}^d} |\nabla u|^p\dx$. By
Example~\ref{ex1}, the functional  $Q_1$ admits a ground state
$\varphi=\mathrm{constant}$ in $\mathbb{R}^d$.

Let $Q_0$ be a functional of the form \eqref{Q} satisfying
$Q_0\geq 0$ on $C_0^\infty(\mathbb{R}^d)$. Let   $\psi\in
W^{1,p}_{\mathrm{loc}}(\mathbb{R}^d)$, $\psi_+\neq 0$ be a
subsolution of the equation $Q_0'(u)=0$ in $\mathbb{R}^d$, such
that $\psi_+\in L^\infty(\mathbb{R}^d)$. It follows from
Theorem~\ref{thm:main} that $\psi$ is the ground state of $Q_0$ in
$\mathbb{R}^d$. In particular, $\psi$ is (up to a multiplicative
constant) the unique positive supersolution and unique bounded
solution of the equation $Q'_0(u)=0$ in $\mathbb{R}^d$. Note that
there is no assumption on the behavior of the potential $V_0$ at
infinity. This result generalizes some striking Liouville theorems
for Schr\"odinger operators on $\mathbb{R}^d$ that hold for
$d=1,2$ and $p=2$ (see \cite[theorems~1.4--1.6]{Pl}).
 }
 \end{example}
\begin{example}\label{ex2}{\em
Let $d>1$, $d\neq p$, and $\Omega:=\mathbb{R}^d\setminus \{0\}$ be
the punctured space. Let $c^*_{p,d}:=|{(p-d)}/{p}|^p$ be the Hardy
constant, and consider the functional
\begin{equation}\label{Hardy_fun}
  Q(u):=\int_\Omega \left(|\nabla u|^p-
  c^*_{p,d}\frac{|u|^p}{|x|^p}\right)\dx
 \qquad u\in \core.
\end{equation}
By  Hardy's inequality, $Q$ is nonnegative on $\core$.  The proof
of Theorem 1.3 in \cite{PS} shows that $Q$ admits a null sequence.
It can be easily checked that the function $v(r):=|r|^{(p-d)/p}$
is a positive solution of the corresponding radial equation:
$$-|v'|^{p-2}\left[(p-1)v''+\frac{d-1}{r}v'\right]-c^*_{p,d}\frac{|v|^{p-2}v}{r^p}
=0\qquad r\in (0,\infty).$$ Therefore, $\varphi(x):=|x|^{(p-d)/p}$
is the ground state of the equation
\begin{equation}\label{eqHardy}
-\Delta_p u- c^*_{p,d}\frac{|u|^{p-2}u}{|x|^{p}}=0 \qquad \mbox{in
} \Omega.
\end{equation}
Note that $\varphi\not \in W^{1,p}_{\mathrm{loc}}(\R^d)$ for
$p\neq d$. In particular, $\varphi$ is not a positive
supersolution of the equation $\Delta_pu=0$ in $\R^d$.
%

Let $Q_0$ be a functional of the form \eqref{Q} satisfying
$Q_0\geq 0$ on $\core$. Let   $\psi\in
W^{1,p}_{\mathrm{loc}}(\Omega)$, $1<p<\infty$, $p\neq d$,
$\psi_+\neq 0$ be a subsolution of the equation $Q_0'(u)=0$ in
$\Omega$, satisfying \be\label{HardyPsi} \psi_+(x)\le
C|x|^{(p-d)/p},\; x\in\Omega.\ee When $p>2$, we require in
addition that the following inequality is satisfied \be
\label{HardyGrad} \psi_+(x)^2|\nabla\psi_+(x)|^{p-2}\le
C|x|^{2-d},\; x\in\Omega. \ee It follows from
Theorem~\ref{thm:main}, Remark~\ref{remgradcond} and
Remark~\ref{remNogradcond} that $\psi$ is the ground state of
$Q_0$ in $\Omega$. The reason why \eqref{HardyGrad} is stated only
for $p>2$ hinges on the fact that for  $p\le 2$, \be
\label{HardyGS} C^{-1} Q_0(\varphi w)\le
E_1^\varphi(w)=\int_{\Omega}|x|^{p-d}|\nabla w|^p\dx, \ee and for
all $p>1$ the functional $E_1^\varphi$  has a ground state
$\mathbf{1}$. The null sequence convergent to this ground state is
given by \cite{PS}, relation (2.2) with $R\to\infty$. Therefore,
Remark~\ref{remNogradcond} applies. }
 \end{example}
Next, we present a family of functionals $Q_0$ for which the conditions of
 Example~\ref{ex2} are satisfied.
\begin{example}
{\em Let $d\ge 2$, $1<p<d$, $\alpha\ge 0$, and
$\Omega:=\mathbb{R}^d\setminus \{0\}$. Let
$$
W_\alpha(x):=-\left(\dfrac{d-p}{p}\right)^p
\dfrac{\dfrac{\alpha\, d\,p}{d-p}+|x|^\frac{p}{p-1}}
{\left(\alpha+|x|^\frac{p}{p-1}\right)^p}\;.
$$
Note that if $\alpha=0$ this is the Hardy potential as in the
Example~\ref{ex2}. If $Q_0$ is the functional \eqref{Q} with
the potential $V:=W_\alpha$, then
$$
\psi_\alpha(x):=\left(\alpha+|x|^\frac{p}{p-1}\right)^{-\frac{(d-p)(p-1)}{p^2}}
$$
is a solution of $Q_0'(u)=0$ in $\Omega$, and therefore $Q_0\ge 0$
on $\core$. Moreover, one can use the calculations of
Example~\ref{ex2} to show that $\psi_\alpha$ is a ground state of
$Q_0$. Indeed, we note first that $\psi=\psi_\alpha$ satisfies
\eqref{HardyPsi}. If $\frac{d}{d-1}<p<d$, then $\psi_\alpha$
satisfies also \eqref{HardyGrad} and therefore, it is a ground
state in this case. In the remaining case $p\le\frac{d}{d-1}\le 2$,
Example~\ref{ex2} concludes that $\psi_\alpha$ is a ground state
from the property of the functional \eqref{HardyGS}. }
\end{example}

\mysection{The simplified energy}\label{secsimpl}

In this section we give examples showing that none of the terms in
the simplified energy \eqref{p>2} for $p>2$ is dominated by the
other, so that \eqref{p>2} cannot be further simplified. In
particular, neither condition  \eqref{ineq:phi-psi} nor condition
\eqref{ineq:gradphi-gradpsi} in Theorem \ref{thm:main} can be
omitted.

Let $p>2$, and fix $v>0$, $v\in C^1_{\mathrm{loc}}(\Omega)$. For
$w\in\core$ we denote \be \label{Lagr-1} E^v_1(w):=\int_{\Omega}
v^p|\nabla w|^p\dx,
 \ee
and \be \label{Lagr-2} E^v_2(w):=\int_{\Omega}v^2|\nabla v|^{p-2}
w^{p-2}|\nabla w|^2\dx. \ee
\begin{remark}
\label{rem:irreducible} {\em Suppose that $V=0$ in $\Omega$, and
assume that $Q$ has a weighted spectral gap in $\Omega$. In this
case, the constant function $v=\mathbf{1}$ is a positive solution
of \eqref{groundstate} which is not a ground state in $\Omega$.
Clearly,  $E_2^{v}= 0$ on $\core$. Therefore, the inequality
$Q(vu)\le CE_2^v(u)$ for $u\in \core$ is generally false, or in other
words, the first term in the simplified energy \eqref{p>2} is
necessary. The above assumptions are satisfied if $Q'$ is the
$p$-Laplacian operator, and either
$\mathrm{int}\,(\R^d\setminus\Omega)\neq\emptyset$, or
$\Omega=\R^d$ and $p< d$ (see Example~\ref{ex1}).
 }
 \end{remark}
In the following proposition we restrict our consideration to the
case $\Omega=\R^d$ and a radial positive solution $v$.
\begin{proposition}\label{prop:irreducible}
 Let $\Omega=\R^d$ and $p>2$.
\begin{itemize}
\item There exists a positive continuous radial function $\varphi$
on $\R^d$ such that the simplified energy $E_1^\varphi(w)+
E_2^\varphi(w)$  has a weighted spectral gap, but there exist a
sequence $\{u_k\}\subset C_0^\infty(\R^d)$ with $u_k\ge 0$, and an
open set $B\Subset\R^d$
 such that
$E_1^\varphi(u_k)\to 0$ and $\int_B|u_k|^p\dx=1$.
\item Moreover,
there exists a positive radial continuous function $\psi$ on
$\R^d$ such that the simplified energy $E_1^\psi(w)+ E_2^\psi(w)$
has a weighted spectral gap, but there exist a sequence
$\{v_k\}\subset C_0^\infty(\R^d)$ with $v_k\ge 0$, and an open set
$B\Subset\R^d$ such that $E_2^\psi(v_k)\to 0$ and
 $\int_B|v_k|^p\dx=1$.
\end{itemize}
\end{proposition}
\begin{proof}
\par\noindent {\bf Step 1.} Let $Q$ be a functional of the form
\eqref{Q} on $\R^d$ with a radial potential $V$. By the proof of
Theorem~\ref{pos} (see, \cite[Theorem~2.3]{ky3}), it follows that
a {\em radial} functional $Q$ is nonnegative on $C_0^\infty(\R^d)$
if and only if the equation $Q'(u)=0$ admits a positive {\em
radial} solution in $\R^d$. On the other hand, from the standard
rearrangement argument it is evident that $Q$ has a weighted
spectral gap if and only if it has a radial weighted spectral gap
with a radial potential $W$. Therefore, $Q$ has a weighted
spectral gap if and only if there exists a positive continuous
radial potential $W$ such that the Euler-Lagrange equation for the
functional $Q(u)-\int_{\R^d}W|u|^p\dx$ has a positive {\em radial}
solution. But such a solution is a (radial) supersolution of the
equation $Q'(u)=0$ in $\R^d$ which is not a solution. Therefore,
it is sufficient to consider the restrictions of $Q$ and $Q'$ to
radial functions. (So, in fact, we are dealing with a
one-dimensional problem.)
\par\noindent {\bf Step 2.} We establish for a radial function $\varphi$ a criterion for the
existence of null sequences for $E^\varphi_1$ and $E^\varphi_2$
using a change of variable. First, for a positive continuous
function $\varphi$ on $[0,\infty)$ define $\rho_1(r)$ by
\be\label{rhor} \rho_1(r) := \int_1^r
\left[\varphi^p(s)s^{d-1}\right] ^{-1/(p-1)}\ds\;.\ee Assume
further that $p\leq d$, then $\rho_1$ is well-defined on
$(0,\infty)$ and $\rho_1(0)=-\infty$. Since
$$\frac{\mathrm{d}r}{\mathrm{d}\rho_1}=(\varphi^pr^{d-1}) ^{1/(p-1)},$$
it follows  that for a radially symmetric $w\in C_0^\infty(\R^d)$
we have \be E^\varphi_1(w)=\int_{\R^d} \varphi^p|\nabla
w|^p\dx=C_d\int_{-\infty}^{M_1}|w'(\rho_1)|^p\drho_1, \ee where
\be\label{1} M_1=M_1^\varphi: =\int_1^\infty
 (\varphi^p r^{d-1}) ^{-1/(p-1)}\dr.\ee
Recall that from Example~\ref{ex1} it follows that for $p>1$ the
$p$-Laplacian on $(a,b)$  admits a ground state if and only if
$(a,b)=\R$. Therefore, $E^\varphi_1$ has a null sequence if and
only if $M_1^\varphi=\infty$ and $p\leq d$ (cf. \cite[Theorem
3.1]{Murata86}).
\par
Consider now the functional $E^\varphi_2$. The substitution
$u:=w^{p/2}$ implies that \be
E^\varphi_2(w)=\int_{\R^d}\varphi^2\,|\nabla \varphi|^{p-2}
w^{p-2}|\nabla w|^2\dx= (2/p)^{2}\int_{\R^d}\varphi^2
|\nabla\varphi|^{p-2}|\nabla u|^2\dx.\ee Let \be\label{rhor1}
\rho_2(r) := \int_1^r \varphi(s)^{-2}|\varphi'(s)|^{2-p}s^{1-d}
\ds \;,\ee and assume further that $d\geq 2$, then $\rho_2$ is
well-defined function on $(0,\infty)$ and $\rho_2(0)=-\infty$.

Using spherical coordinates for radial $w$, and then the
substitution $\rho_2$, we obtain
 $$\int_{\R^d}\varphi^2\,|\nabla \varphi|^{p-2}
w^{p-2}|\nabla w|^2\dx=C_d\int_{-\infty}^{M_2}|
u'(\rho_2)|^2\drho_2,$$ where \be\label{2}
M_2=M_2^\varphi=\int_1^\infty
 \varphi^{-2}|\varphi'|^{2-p}r^{1-d}\dr.
\ee Therefore, $E^\varphi_2$ has a null sequence (for $p>2$ and
$d\geq 2$) if and only if $M_2^\varphi=\infty$ (cf. \cite[Theorem
3.1]{Murata86}).

Therefore, in order to prove the proposition it is sufficient to
find two positive radial functions $\varphi$ and $\psi$ satisfying
$M_1^\varphi=\infty$ and $M_2^\varphi<\infty $, while
$M_1^\psi<\infty $ and $M_2^\psi=\infty$.

\vskip 4mm

\par\noindent {\bf Step 3.} Let us simplify now \eqref{1} and \eqref{2}
in order to investigate  when $M_j^\varphi$ are finite or infinite
for a specific $\varphi$ and $j=1,2$. Without loss of generality
we assume that the integration in \eqref{1}  and \eqref{2} is from
$r_0$ to $\infty$, where $r_0\gg 1$. We set first
$\varphi(r):=r^{1-d/p}\eta(r)$, where $2<p<d$. Then \eqref{1}
becomes, up to a constant multiple, \be M_1=\int_{r_0}^\infty
 (\varphi^pr^{d-1}) ^{-1/(p-1)}\dr=\int_{{r_0}}^\infty
\eta^{-p/(p-1)}r^{-1}\dr.\ee
 Set now
$$\eta(r):=\left[t^{(p-1)/(p-2)}(\log t)^\gamma\right]^{(p-2)/p},
\quad \mbox{ where } t:=\log r, \mbox{ and } \gamma >0.$$   Then
we have
\begin{multline*}\label{1a} M_1=\int_{r_0}^\infty \eta^{-p/(p-1)}r^{-1}\dr=
\int_{t_0}^\infty \left[t^{(p-1)/(p-2)}(\log
t)^\gamma\right]^{(2-p)/(p-1)}\dt=\\\int_{t_0}^\infty \frac{(\log
t)^{\gamma(2-p)/(p-1)}}{t}\dt . \end{multline*} On the other hand
\be\label{2a12} M_2=\int_{r_0}^\infty
 \varphi^{-2}|\varphi'|^{2-p}r^{1-d}\dr=\int_{r_0}^\infty \left|\frac{p-d}{p}+
\frac{r\eta'}{\eta}\right|^{2-p}\eta^{-p}r^{-1}\dr. \ee

Denote \be\label{2b1} \tilde{M}_2:=\int_{r_0}^\infty
\eta^{-p}r^{-1}\dr=\int_{t_0}^\infty t^{1-p}(\log
t)^{\gamma(2-p)}\dt. \ee Since  $r|\eta'|/\eta\ll 1$, it follows
from  \eqref{2a12} that there exist $C>0$ such that
 \be \label{m2}
C^{-1}\tilde{M}_2\leq M_2\leq C\tilde{M}_2.
 \ee
Consequently, for $0<\gamma\leq (p-1)/(p-2)$ and $2<p<d$, we have
$M_1^\varphi=\infty$ and $M_2^\varphi<\infty$.

\vskip 3mm

On the other hand, fix  $\beta\in \mathbb{R}$, $\beta\neq 0$,  and
let $\psi:(1,\infty)\to [1,\infty)$,  be a smooth monotone
function such that $\psi(r)\asymp r^\beta$, and such that $\psi'$
satisfies for any $n=1,2,\dots,$
$$|\psi'(r)|= \left\{
\begin{array}{ll}
    \mathrm{e}^{-r} & r\in [2n+1/4, 2n+3/4], \\[2mm]
    |\beta|r^{\beta-1} & r\in [2n+1, 2n+2].
\end{array}%
\right.$$
Therefore, if $\beta>(p-d)/p$, then $M_1^\psi<\infty$.

Consider now $M_2^\psi$, and recall that $2<p\leq d$.
Consequently, there exist $\varepsilon>0$ and $C>0$ such that the
integrand of $M_2^\psi$ satisfies
$$\psi^{-2}|\psi'|^{2-p}r^{1-d}\geq C\mathrm{e}^{\varepsilon
r}$$ on a set of infinite measure. Hence, $M_2^\psi=\infty$.
\end{proof}
\begin{remark}\label{differentiation}{\em
The proof above takes into account that
\eqref{1}, \eqref{2} both converge or both diverge when
$\varphi(r)$ is of the form $r^{\alpha}(\log r)^\beta$, and the
differentiation occurs only with respect to $\gamma$ for
$\varphi(r)=r^{1-d/p}(\log r)^{(p-1)/p}(\log\log r)^\gamma$.
 }
 \end{remark}

\mysection{Application: ground state of the linearized
functional}\label{secappl} We consider the linearized problem
associated with the functional $Q\ge 0$. Let $\varphi$ be a
positive solution of the equation $Q'(u)=0$ in $\Omega$, and let
\be
\label{linearized} \mathbf{a}[u]:=\int_\Omega \left(|\nabla
\varphi|^{p-2}|\nabla u|^2+V(x)\varphi^{p-2}u^2\right)\dx.\ee
\begin{proposition} Let $\varphi$ be a positive solution of the
equation $Q'(u)=0$ in $\Omega$ satisfying $\nabla \varphi\neq 0$.
\begin{enumerate}
\item If $p>2$ and $\varphi$ is a ground state of $Q$, then $\varphi$ is a ground state of
 $\mathbf{a}$.
\item  If $p<2$ and $\varphi$ is a ground state of $\mathbf{a}$,
then $\varphi$ is a ground state of $Q$.
\end{enumerate}
\end{proposition}
\begin{proof} Consider first the case $p>2$. Assume that
$\varphi$ is a ground state of $Q$. Let $\{u_k\}$ be a null
sequence of nonnegative functions, and let $w_k:=u_k/\varphi$. By
\eqref{p>2}, \be \label{nullQ}
\int_\Omega\varphi^2|\nabla\varphi|^{p-2} w_k^{p-2}|\nabla w_k|^2
\dx\to 0. \ee
 Set $v_k:=w_k^{p/2}$. Then \be \label{nullH}
\int_\Omega\varphi^2|\nabla\varphi|^{p-2}|\nabla v_k|^2 \dx\to 0
\ee which by \eqref{e2} yields
$$\mathbf{a}[\varphi v_k]\to 0.$$
Taking into account Remark~\ref{remc1}, we conclude that
$\{\varphi v_k\}$ is a null sequence for $\mathbf{a}$.
\par
The case $p<2$ is similar. If $\{z_k\}$ is a null sequence  of
nonnegative functions for the form $\mathbf{a}$, then
\eqref{nullH} is satisfied with $v_k:=z_k/\varphi$. This implies
\eqref{nullQ} with $w_k=v_k^{2/p}$,  which by \eqref{p<2} yields
$Q(u_k)\to 0$ with $u_k=\varphi w_k$. Therefore, $\{u_k\}$ is a
null sequence for $Q$ and the proposition is proved.
\end{proof}
\begin{center}
{\bf Acknowledgments} \end{center} Parts of this research were
done while K.~T. was visiting the Technion, Y.~P. was visiting
University of Crete and Uppsala University and A.T. was visiting
Uppsala University. Y.~P. and A.~T. acknowledge partial  support
by the RTN European network Fronts--Singularities,
HPRN-CT-2002-00274. The visit to Technion was supported by the
Fund for the Promotion of Research at the Technion and the visits
to Uppsala were supported by the Swedish Research Council.
%


\begin{thebibliography}{00} {\frenchspacing
%
\bibitem{Ag} S.~Agmon,  Bounds on exponential decay of eigenfunctions of
Schr\"odinger operators, {\em in}  ``Schr\"odinger Operators"
(Como, 1984), pp. 1--38, Lecture Notes in Math., 1159, Springer,
Berlin, 1985.

\bibitem{AH1} W.~Allegretto, and Y.~X.~Huang, A Picone's identity for the
$p$-Laplacian and applications, {\em Nonlinear Anal.} {\bf 32}
(1998), 819--830.
%
\bibitem{AH2} W.~Allegretto, and Y.~X.~Huang, Principal eigenvalues and Sturm
comparison via Picone's identity, {\em  J. Differential Equations}
{\bf 156} (1999), 427--438.
%
%
\bibitem{BFT} G.~Barbatis, S.~Filippas, and A.~Tertikas,
A unified approach to improved $L^p$ Hardy inequalities with best
constants, {\em  Trans. Amer. Math. Soc.} {\bf 356} (2004),
2169--2196.
%
\bibitem{DS} J.~I.~Diaz, J.~E.~Sa\'{a},  Existence et unicit\'{e} de
solutions positives pour certaines \'{e}quations elliptiques
quasilin\'{e}aires, {\em C. R. Acad. Sci. Paris Ser. I Math.} {\bf
305} (1987), 521--524.

%
%
%
%
%
%
\bibitem{HKM} J.~Heinonen, T.~Kilpel\"{a}inen, and O.~Martio,
``Nonlinear Potential Theory of Degenerate Elliptic Equations",
Oxford Mathematical Monographs, Oxford University Press, New York,
1993.
%
\bibitem{MP} Mitidieri, and S.~I.~Pokhozhaev, Some generalizations of
Bernstein's theorem, {\em Differ. Uravn.} {\bf 38} (2002),
373--378; translation in {\em Differ. Equ.} {\bf 38} (2002),
392--397.

\bibitem{Murata86} M.~Murata, Structure of positive solutions to $(-\Delta
+V)u=0$ in $\mathbb{R}^n$, {\em Duke Math. J.} {\bf 53} (1986),
869--943.
%
\bibitem{Pl} Y.~Pinchover, A Liouville-type theorem for Schr\"odinger
operators, to appear in Comm. Math. Phys.,\\
http://arxiv.org/PS\_cache/math/pdf/0511/0511039.pdf
%
\bibitem{ky3} Y.~Pinchover, and K.~Tintarev, Ground state
alternative for $p$-Laplacian with potential term, to appear in
Calc. Var. Partial Differential Equations {\bf 28} (2207), 179--201.
%
\bibitem{PS} A.~Poliakovsky, and I.~Shafrir, Uniqueness of
positive solutions for singular problems involving the
$p$-Laplacian, {\em Proc. Amer. Math. Soc.} {\bf 133} (2005),
2549--2557.
%
\bibitem{Serrin1} J.~Serrin, Local behavior of solutions of quasi-linear
equations, {\em Acta Math.} {\bf 111} (1964), 247--302.
%
\bibitem{Serrin2} J.~Serrin, Isolated singularities of solutions of quasi-linear
equations, {\em Acta Math.} {\bf 113} (1965), 219--240.
%
\bibitem{Shafrir} I.~Shafrir, Asymptotic behaviour of
minimizing sequences for Hardy's inequality, {\em  Commun.
Contemp. Math.} {\bf 2} (2000),  151--189.
%
\bibitem{T} P.~Tolksdorf,  Regularity for a more general class of
quasilinear elliptic equations, {\em J. Differential Equations}
{\bf  51} (1984), 126--150.
%
 }\end{thebibliography}
\end{document}